\documentclass[a4paper,fleqn]{cas-dc}

\usepackage{caption}
\usepackage[ruled,linesnumbered]{algorithm2e}
\usepackage{float}
\usepackage{algpseudocode}
\usepackage{amsmath}
\usepackage{subfigure}
\usepackage{graphicx}   
\usepackage{subcaption}
\usepackage{cite}
\usepackage{float}
\usepackage[square]{natbib}
\usepackage{svg}

\setcitestyle{numbers,square}
\def\tsc#1{\csdef{#1}{\textsc{\lowercase{#1}}\xspace}}
\tsc{WGM}
\tsc{QE}
\tsc{EP}
\tsc{PMS}
\tsc{BEC}
\tsc{DE}

\begin{document}
\let\WriteBookmarks\relax
\def\floatpagepagefraction{1}
\def\textpagefraction{.001}

\renewcommand\thesubsubsection{\Alph{subsubsection}}

\shorttitle{Leveraging social media news}

\shortauthors{Hongyaoxing Gu et~al.}

\title [mode = title]{A method for accelerating low precision operations by sparse matrix multiplication}                      
\tnotemark[1,2]



%
\author[1,2]{Hongyaoxing Gu}[type=editor,
                        auid=000,bioid=1,
                        orcid=0009-0007-9276-9536]

\cormark[1]

\fnmark[1]

\ead{771582678@qq.com}


\credit{Conceptualization of this study, Methodology, Software}

\affiliation[1]{organization={Institute of Software Chinese Academy of Sciences},
    city={BeiJing},
    country={China}}




\credit{Data curation, Writing - Original draft preparation}

\affiliation[2]{organization={University of Chinese Academy of Sciences},
    city={Bei Jing},
    country={China}}



\cortext[cor1]{Corresponding author}
\cortext[cor2]{Principal corresponding author}



\begin{abstract}
In recent years, the fervent demand for computational power across various domains has prompted hardware manufacturers to introduce specialized computing hardware aimed at enhancing computational capabilities. Particularly, the utilization of tensor hardware supporting low precision has gained increasing prominence in scientific research. However, the use of low-precision tensor hardware for computational acceleration often introduces errors, posing a fundamental challenge of simultaneously achieving effective acceleration while maintaining computational accuracy.

This paper proposes improvements in the methodology by incorporating low-precision quantization and employing a residual matrix for error correction and combines vector-wise quantization method.. The key innovation lies in the use of sparse matrices instead of dense matrices when compensating for errors with a residual matrix. By focusing solely on values that may significantly impact relative errors under a specified threshold, this approach aims to control quantization errors while reducing computational complexity. Experimental results demonstrate that this method can effectively control the quantization error while maintaining high acceleration effect.The improved algorithm on the CPU can achieve up to 15\% accuracy improvement while 1.46 times speed improvement.
\end{abstract}


\begin{highlights}
\item A hybrid precision quantization method
\item Sparsity is introduced into low precision to improve accuracy
\end{highlights}


\maketitle

\section{Introduction}
In the field of deep learning, neural networks typically employ 32-bit or 64-bit floating-point numbers to represent weights and activation values, which incurs substantial memory usage and high computational costs. Low-precision quantization is a computational optimization technique which can apply to deep learning, aiming at diminishing model storage requirements and computational expenses while preserving model performance.

Tensors, as a high-dimensional generalization of matrices, have become a pivotal data structure in intelligent applications such as deep learning \cite{ref1}.  As deep neural networks are increasingly deployed, various types of tensor-specific hardware have been introduced to enhance the performance and energy efficiency of tensor computations.  Prominent among these heterogeneous platforms are Google's Tensor Processing Unit (TPU) \cite{ref2}, Intel's Neural Network Processor (NNP), Neural Processing Units (NPU) \cite{ref3, ref4}, and NVIDIA GPU \cite{ref5}.  These computing devices incorporate dedicated tensor units known as tensor cores\cite{ref6}.  For instance, the Nvidia A100 introduces a potent third-generation tensor core, exhibiting significantly higher throughput than the V100.  It also features comprehensive support for DL and HPC data types, along with new sparsity capabilities, resulting in doubled throughput.

Nvidia's architectures, including Tesla and Ampere, feature tensor cores that support low-precision numerical computations.  This involves utilizing different precision values for input and output in mixed-precision calculations.  Computational operations are performed using low precision, with the final results stored in high precision.  This approach enhances speed through low-precision calculations while ensuring that the results do not suffer from numerical overflow.

For FP16/FP32 mixed-precision DL, A100 Tensor Core performance is 2.5 times that of V100, increasing to 5 times with added sparsity.  Furthermore, the acceleration of INT8, INT4, and binary Tensor Core implementations supports deep learning inference.  Table 1 shows the speed support for different low-precision tensor core calculations across Nvidia architectures, as indicated by official Nvidia evaluations\cite{ref7, ref8, ref9}.

\begin{table}[htbp]
  \centering
  \caption{Low-precision computing types and performance supported by tensorcore on different GPU computing platforms}
    \begin{tabular}{|c|c|c|}
    \toprule
    Arch and Model & Compute-Accumulator & Performance \\
    \midrule
    Volta V100 & FP16-FP32 & 125 TFLOPS \\
    \midrule
    \multirow{3}[6]{*}{Tesla T4} & FP16-FP32 & 65  TFLOPS \\
\cmidrule{2-3}          & INT8-INT32 & 130 TOPS \\
\cmidrule{2-3}          & INT4-INT32 & 260 TOPS \\
    \midrule
    \multirow{4}[8]{*}{Ampare A100} & FP16-FP32 & 312 TFLOPS \\
\cmidrule{2-3}          & INT8-INT32 & 624 TOPS \\
\cmidrule{2-3}          & INT4-INT32 & 1248 TOPS \\
\cmidrule{2-3}          & INT1-INT32 & 4992 TOPS \\
    \bottomrule
    \end{tabular}%
  \label{tab:addlabel}%
\end{table}%

Besides, CPUS are starting to support more high performance low-precision operations. Intel ® Xeon ® scalable processors are capable of using low-precision arithmetic operations (INT8) in 8-bit integer format. The method achieves lower latency than 32-bit single-precision floating-point arithmetic (FP32) workloads while maintaining high accuracy. The built-in Intel ® Advanced Vector Extension 512 (Intel ® AVX-512) and Vector Neural Network Instructions (VNNI) further improve INT8 inference performance \cite{ref43}.

In the next generation of tensor core computing devices, low-precision operations based on tensors have gained widespread support.  Tensor-specific hardware enables the execution of tensor computations in a single operation, significantly boosting the throughput and efficiency of linear algebra operations, such as matrix multiplication.  To leverage the efficient computational capabilities of low precision, quantization methods are employed \cite{ref10,ref28}.  Low-precision quantization reduces the number of bits required to represent weights and activation values, thereby substantially decreasing model storage requirements and computational costs.  Typically, low-precision quantization involves converting floating-point representations to fixed-point or integer representations, such as using 8-bit integers or fewer bits.  While this conversion introduces some precision loss, in many cases, this loss can be compensated for through various techniques to maintain model performance.

Low-precision quantization is particularly useful in resource-constrained environments such as hardware accelerators, edge devices, and mobile devices, as it significantly reduces the storage and computational overhead of models, making them more suitable for deployment on these devices.  Beyond the field of machine learning, there is a growing recognition that low precision can also accelerate traditional high-performance applications \cite{ref11}.

This paper revolves around the following question: How can we take advantage of the high computational performance of low precision while ensuring that the error of the result is acceptable? 
Based on this problem, the contributions made in this paper are as follows:
\begin{itemize}
    \setlength{\itemsep}{0pt} 
    \item The residual error repair method is analyzed, and a method to control the calculation amount of low-precision matrix multiplication by threshold value is given.
    \item Sparse matrix is introduced in the calculation process, and the calculation amount is reduced by sparse matrix multiplication within the range of error threshold.
    \item Designed a high-performance low-precision quantization algorithm by Cutlass and Cusparse on GPU, and proved the effectiveness of the algorithm through a series of experiments on Nvidia-A100.
\end{itemize}

\section{Motivation}

Leveraging tensor hardware for optimizing tensor computations will be a key factor in enhancing application computational performance. In linear algebra, matrix multiplication is a common operation, and tensors are widely used to represent multidimensional data. In deep learning, the majority of computations arise from tensor operations, with tensor multiplication constituting the majority of these computations (convolution operations can be transformed into matrix multiplication through techniques like im2col\cite{ref12}). There exists a close relationship between tensor computations and matrix multiplication, where matrix multiplication is a special case of tensor multiplication. In matrix multiplication, the result of multiplying two matrices is a new matrix. Similarly, in tensor multiplication, applicable to higher-order tensors, the result of multiplying two tensors is a new tensor, with elements representing combinations of elements from the original tensors. Therefore, optimizing matrix multiplication through tensor cores can accelerate tensor computation applications.

In Section 2.1, we will explore the performance of low-precision operations on advanced computing devices and the acceleration effects brought about by sparse matrices. In Section 2.2, we will delve into the process of low-precision quantization and dequantization for matrices, demonstrating the use of quantization in the low-precision execution of matrix multiplication $C=A * B $. Then using cuBlas\cite{ref23} and CuSparse\cite{ref25} libraries as examples, we will present a performance comparison of low-precision matrix multiplication. In Section 2.3, we will introduce existing research on low-precision quantized matrix multiplication, introducing its advantages and limitations.

\subsection{Low precision computation and sparse matrix acceleration}
\subsubsection{Low precision GEMM}
Due to the high computational performance of low-precision calculations, current computing acceleration hardware is gradually supporting a wider range of low-precision data types for tensor computations.   The A100, for instance, introduces new low-precision calculation types, including int8, int4, int1(equivalent to bool).   Additionally, Tensor Cores on NVIDIA GPUs support low-precision-mixed precision computation.   This method involves using different precisions for input and output, typically producing high-precision output to prevent potential data overflow issues arising from low-precision calculations.

\begin{figure}
	\centering
\includegraphics[width=1\linewidth]{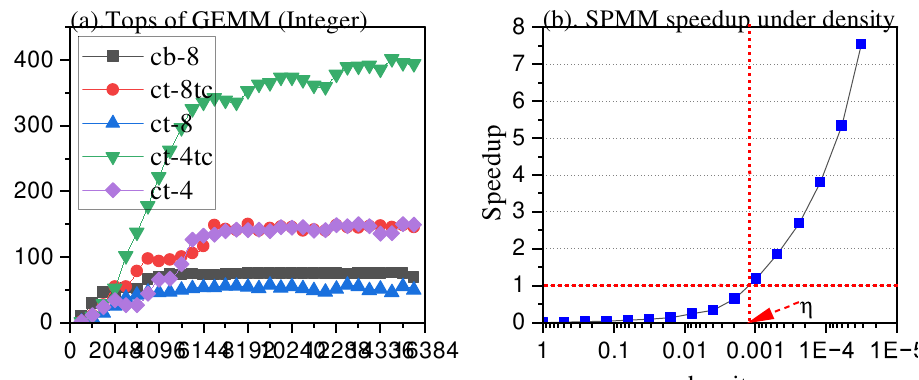}
\includegraphics[width=1\linewidth]{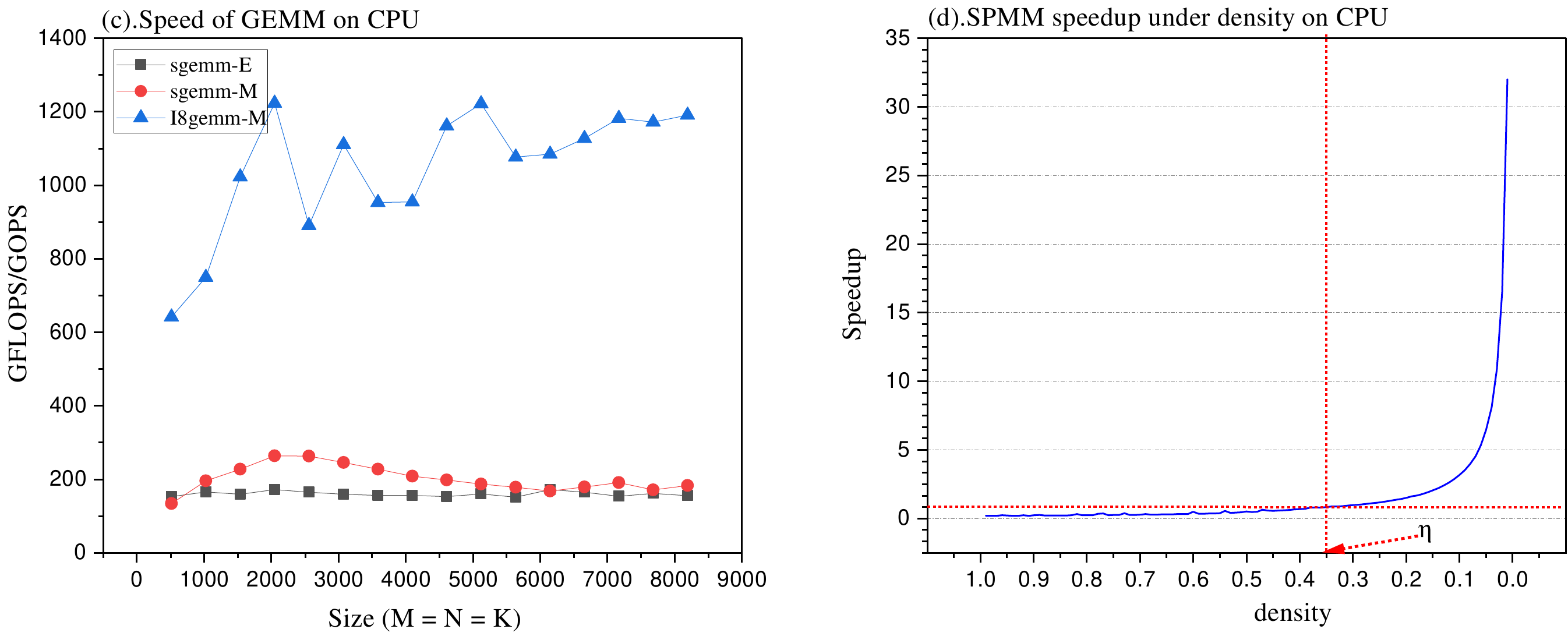}
	\caption{Figures(a) includes the efficiency of  (int8-int32), and (int4-int32) executions.   Here, 'cb' represents Cublas, 'ct' represents Cutlass, and 'tc' denotes the use of Tensor Cores. Figures(c) show the execution time of GEMM with different precision in different math libraries under Intel(R) Xeon(R) Platinum 8163. Where E stands for EIGEN mathematical library and M stands for MKL mathematical library. I8 stands for int8 as the calculation accuracy and int32 as the result accuracy(Eigen does not support the calculation).  Figures(b) and Figures(d) shows the acceleration of sparse matrix multiplication (SPMM) over dense matrix multiplication (GEMM) on GPU/CPU. $\eta$ indicates the density of SPMM with the same speed as GEMM. 
    }
	\label{FIG:diffper}
\end{figure}

Cublas and Cutlass\cite{ref23, ref24} are mainstream mathematical libraries for performing general matrix-multiplication (GEMM) calculations on Nvidia GPUs.   Cutlass, in comparison to Cublas, supports a broader range of low-precision data types and matrix operator designs, including int1, int4, int8, float16, and bfloat16. This paper focuses on the acceleration of matrix calculation by low precision shaping. The Fig.\ref{FIG:diffper}(a) below shows the execution speeds of GEMM (Matrix multiplication) with different precision for integers numbers on the Nvidia A100-40GB.

The Intel MKL Math Core Library\cite{ref44} is a highly optimized and extensively threaded set of math routines that currently support integer matrix multiplication down to int8. And Eigen\cite{ref45}  is a high-level C ++ library that effectively supports linear algebra, matrix and vector operations, numerical analysis and their related algorithms. \ref{FIG:diffper}(c) shows the execution speeds of GEMM (Matrix multiplication) with different precision for integers numbers on the Intel(R) Xeon(R) Platinum 8163 CPU.  The low-precision int8 matrix multiplications shown in the figure can achieve up to 8 times faster than single-precision floating-point matrix multiplications

In the presented results, it is evident that the computational efficiency of GEMM operations at low precision is significantly enhanced  in both GPU and CPU.


\subsubsection{Sparse matrix acceleration}
In  matrix computations, employing sparse matrices for matrices predominantly composed of zero elements proves to be a more efficient strategy.  Sparse matrices, in contrast to dense matrices where the majority of elements are non-zero, offer notable advantages in terms of storage and computation.  The presence of numerous zero elements within sparse matrices allows for skipping these elements during computations, thereby reducing computational complexity.  This attribute results in significant performance improvements, particularly in numerical computations and linear algebra operations such as matrix multiplication.  Sparse matrices demonstrate enhanced efficiency and reduced resource consumption when handling large-scale, high-dimensional data characterized by evident sparsity.  Consequently, they find widespread application in various scientific, engineering, and computational domains.

Taking int8-int32 GEMM (dense matrix multiplication) as a benchmark, the Fig.3 below shows the acceleration ratios of SPMM (sparse matrix multiplication) at varying levels of sparsity.The matrices entered for the above two experiments are randomly generated matrices using the Curand library.

From the experimental results, it is evident that, with sufficiently high sparsity, sparse matrix multiplication can achieve notably superior acceleration ratios.

\subsection{Quantization and Dequantization}
In accordance with the definition provided by the IEEE 754 standard\cite{ref13}, floating-point numbers consist of three components: the sign bit, the exponent bits, and the fraction bits. For an FP32, or single-precision floating-point number, this entails 1 sign bit, 8 exponent bits, and 23 fraction bits. Similarly, an FP16, or half-precision floating-point number, is comprised of 1 sign bit, 5 exponent bits, and 10 fraction bits. Regarding the quantization of FP16, given that both are represented by floating-point numbers, a direct conversion from high precision to low precision floating-point representation is feasible——($a_{fp16}=TypeCast(a_{fp32},Float16)$). 

However, quantization for integer types is not as straightforward. Due to the inherent structure of floating-point numbers, they are not uniformly distributed along the number line but rather denser near zero. In contrast, integers are uniformly distributed along the number line and typically represent a much smaller range than floating-point numbers. Directly rounding floating-point numbers and representing them as integers is not a reasonable approach, as it would result in the loss of all values near the zero point. Therefore, it is necessary to apply mathematical transformations to floating-point numbers before converting them to integers. This can be achieved, for instance, through clustering\cite{ref14} and scaling\cite{ref15,ref29} to quantize floating-point values into integer values. Alternatively, methods such as using KL divergence\cite{ref16, ref17} involve truncating a portion of the original information before mapping, creating a symmetric and well-distributed truncated information, which is then mapped to the integer domain.

Considering a input matrix, it can be abstracted as a one-dimensional array of numbers arranged from smallest to largest along the number line. The goal of quantization is to map the numerical values from the floating-point domain onto the number line in the integer domain. Different quantization methods offer various advantages, but for matrix multiplication, larger absolute values in the original matrix have a more significant impact on the resulting matrix, and errors in vector multiplication accumulate in a specific element of result matrix. If the method of quantization involves truncating larger values using KL divergence, it may introduce considerable errors into the resulting matrix.Moreover, due to the nature of matrix multiplication, results generated from clustering quantization methods\cite{ref27} utilizing exponential functions cannot be easily restored to the original precision. 

We abbreviated quantization as $Q(A)$ and dequantization as $\widetilde{Q(A)}$. To leverage the efficiency of low precision while facilitating the restoration of result precision, we employ a symmetric linear quantization function. Assume the representation with N bits of integers is used,$a_{max}$ is the largest number of matrix A. 

\begin{equation}
\begin{split}
a_{int} = Q(a_{fp},\lambda) = TypeCast(\lambda*a_{fp},intN)
\label{Quantization}
\end{split}
\end{equation}

This function uses the absolute value of the maximum as a threshold, directly mapping this range proportionally to the range of positive and negative integers. Where
\begin{equation}
\label{getlambda}
\begin{split}
 \lambda=(2^{N}-1)/a_{max}
\end{split}
\end{equation}

\begin{figure}
	\centering
\includegraphics[width=1\linewidth]{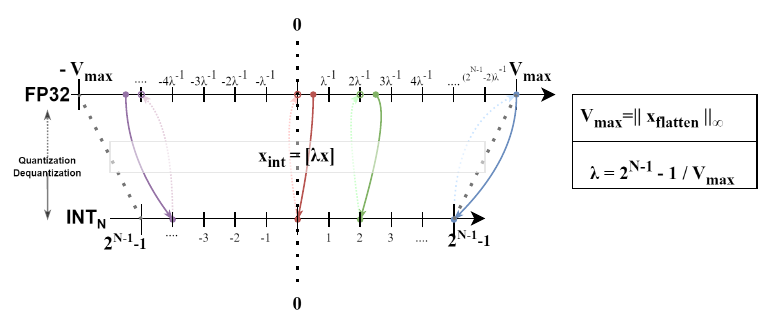}    
	\caption{Representation of integer quantization}
	\label{FIG:7}
\end{figure}

As for the process of dequantization, which is the inverse of quantization, it involves restoring the results represented in integers back to floating-point numbers:
\begin{equation}
\begin{split}
a_{fp} = \widetilde{Q(a_{int},\lambda)} = TypeCast(a_{int}/\lambda,floatN)
\label{Dequantization}
\end{split}
\end{equation}

\subsection{Quantized GEMM and residual-based refinement method}
For GEMM(general matrix multiplication):$D = \alpha AB + \beta C$. For the main part of the calculation, denote$M=AB$. To complete matrix multiplication with low precision, it can be expressed as
\begin{eqnarray}    \label{eq}
    M_{fp32} &=& A_{fp32}*B_{fp32} 
      = \frac{A_{int}}{\lambda_A}*\frac{B_{int}}{\lambda_B}
      = \frac{A_{int}*B_{int}}{\lambda_M} \nonumber \\
      ~&=& \widetilde{Q(A_{int}*B_{int},\lambda_M)}
\end{eqnarray}

Where $\lambda_M = \lambda_A*\lambda_B$. Therefore, we can obtain the original precision matrix by one inverse quantization operation of the low precision M matrix

\begin{figure}
    \centering
    \includegraphics[width=1\linewidth]{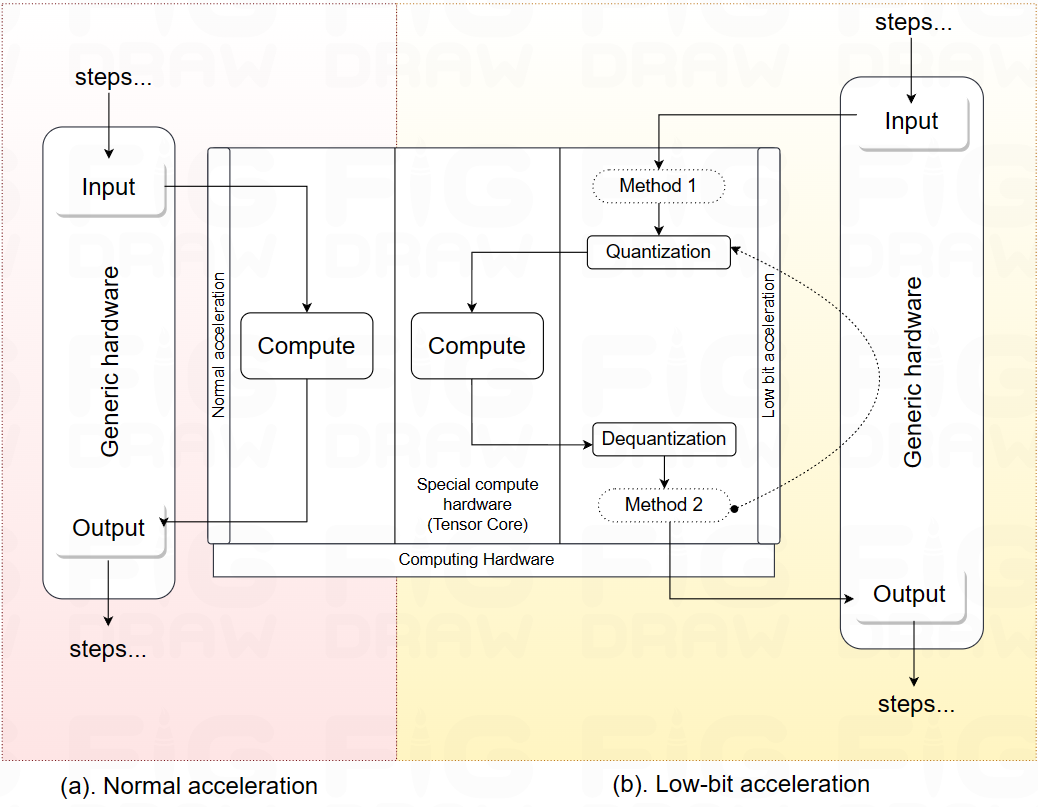}
    \caption{Normal precision acceleration (a), and low precision quantization acceleration processes(b)}
    \label{Fig:5}
\end{figure}

In modern applications (eg. machine learning, numerical computing), the input and output are generally normal precision, and the calculation intensive modules in the program (such as matrix multiplication, convolution operations, etc.) generally account for more than 90\% of the total time of the program. Therefore, the acceleration of compute-intensive modules with low precision can achieve good results. Fig.\ref{Fig:5} shows the flow of computation acceleration using common method and low-precision quantization method.

Compared to the normal precision computing strategy of direct acceleration using computing hardware. Low-precision quantization acceleration first requires quantizing and de-quantizing the inputs, and then using special computing hardware such as tensor core for more efficient computation acceleration. 

Low-precision quantization methods can achieve significant speedup, but the rounding errors introduced by low precision often have a substantial impact on the computation results. Some recent research has employed various techniques to enhance low-precision quantization precision(Fig.\ref{Fig:5}  Method 1, Method 2).

The improvements in low-precision matrix multiplication can be broadly categorized into the following two approaches:
\paragraph{Quantization methods of the input matrix}
The first approach is to use a different method (Fig.\ref{Fig:5} \textit{Method 1}) to change the quantization method before the Count(computationally intensive operation) to improve the accuracy or reduce the calculation time, including the following methods.
\begin{enumerate}
    \item Row and Column wise Quantization\cite{ref18}: This method involves applying different scaling factors $\lambda$ to each row of matrix A and each column of matrix B. The goal is to reduce the variance of the input data, minimizing the difference between the maximum and minimum values and thereby enhancing quantization precision.
    \item Improved Bit Quantization\cite{ref19, ref20}: This approach leverages the bitwise matrix operations supported by the latest GPUs. It enables quantization of inputs with varying bit precision, adapting to different input matrices to enhance quantization precision.
    \item  Matrix Reordering\cite{ref21}: This method involves rearranging the input matrix, expanding matrix B along the rows to transform matrix multiplication into inner products followed by summation. This enables the use of instruction-level algorithms supporting quantization under low precision after the data has been reordered.
\end{enumerate}
The transformation of input matrices can improve quantization precision to some extent, and it does not introduce a significant increase in computational complexity compared to direct quantization. However, it often requires the pre-determination of the quantization bit-width in advance to decide on the computational load, with further fine-tuning based on different inputs. Alternatively, it may only be effective for specific matrix specifications and may perform poorly on matrices with particular characteristics, especially those with large singular values.

\paragraph{Error compensating in output matrix}
The second method (Fig.\ref{Fig:5} \textit{Method 2}) is to compensate the error after the end of the first low-precision calculation, which often requires additional low-precision calculation steps to obtain a considerable accuracy improvement compared with the first method.
The exploration on the improvement of low-precision matrix multiplication mainly focuses on different quantization methods for input matrix(\textit{Method 1}), and there are not many researches on the error compensation of quantized matrix (\textit{Method 2}).   At present, the residual matrix method is mainly used for error compensation\cite{ref22}, and the following is the introduction of this method.

In the aforementioned matrix multiplication quantization, due to $\lambda=(2^{N-1}-1)/a_{max}$ being a floating-point number generated through floating-point computations, this step introduces only minor high-precision floating-point calculation errors which is negligible. Therefore, the quantization operation's error originates solely from the rounding operation in Equation(\ref{Quantization}) , which involves mapping a floating-point number to the integer axis with an offset from its original value.

On the other hand, dequantization, which is the inverse operation of quantization in Equation(\ref{Dequantization}), does not involve rounding when converting an integer to a floating-point number. Therefore, the error incurred in this process is solely due to the inherent error of the floating-point number and can be considered negligible. It can be asserted that dequantization is accurate, and the error exists only in the forward quantization process.

Here's an example:
V = [1, 2.5, 4]. If we want to quantify it using Int8:

Then $\lambda = (2^7-1)/maxV = 31.75$. To number 2.5:\[2.5^{Int} = ToInt(\lambda *a^{float}) =\lfloor31.75*2.5\rfloor= 79 \]

And then we do the dequantization of 79:
$$a^{Float}= ToFloat(\frac{79}{31.75}) = 2.4882$$

Then the error of $Quant(2.5,int4)$ is $2.5-2.4882 = 0.0118$. 
Therefore, for the 8-bit quantization of 2.5, the error is given by the discrepancy caused by $\lfloor79.375\rfloor$.

Performing error calculations for other numbers yields results of $[0.0236, 0.0118, 0]$. We can conceptualize the quantization operation as adding a difference to the original floating-point number, aligning it with the nearest representation on the integer axis. The resulting vector represents the residual vector after int8 quantization of the original vector.

For the low-precision quantization of matrix-matrix multiplication $A\cdot B=C$, we can employ a similar procedure, ultimately obtaining the residual matrix resulting from low-precision quantization. Therefore, the quantization operation for matrix A can be expressed as 
\begin{equation}    \label{eq}
    A^{Fp}=A^{Fp'}+RA^{Fp}
\end{equation}
Where the $A^{Fp'}$ matrix represents the floating-point matrix resulting from the reverse quantization of the integer matrix after quantization, with an offset of $RA^{Fp}$ from the original matrix values.

Applying the same procedure to matrix B, the matrix operation $A\cdot B=C$ can be represented as:

\begin{equation}    \label{eq}
    A^{Fp}\cdot B^{Fp} =(A^{'Fp}+RA^{Fp})*(B^{'Fp}+RB^{Fp})
\end{equation}

The product of two floating-point matrices can be expressed as the multiplication of the quantized integer matrices by the scaling values $\lambda$ of the two matrices.

$$a^{Fp}_{ij}*b^{Fp}_{ij}=\frac{a^{Int}_{ij}}{\lambda_{a}} * \frac{b^{Int}_{ij}}{\lambda_{b}}$$

Therefore, the original matrix multiplication can be obtained by summing up the results of four matrix multiplications:

\begin{equation}
\begin{split}
    A^{Fp}*B^{Fp} &=\frac{A^{Int}\cdot B^{Int}}{\lambda_{a}*\lambda_{b}}+\frac{A^{Int}\cdot R_B^{Int}}{\lambda_{a}*\lambda_{Rb}}  \\
    &+\frac{R_A^{Int}\cdot B^{Int}}{\lambda_{Ra}*\lambda_{b}} 
      +\frac{R_A^{Int}\cdot R_B^{Int}}{\lambda_{Ra}*\lambda_{Rb}} = C^{Fp}
\end{split}
\label{gemm_r}
\end{equation}
This is also a way of precision compensation by using the residual matrix directly. Note that the last term, due to the multiplication of residual matrices, is numerically small compared to the first three terms and can be neglected under certain conditions. We sum up the first three terms to form the low-precision quantized matrix multiplication with residual correction.

Through this method, the error of direct quantized matrix multiplication can be significantly reduced. However, this introduces two new matrix multiplication operations. If the acceleration ratio of direct low-precision quantization is 9, then the theoretical upper limit of the acceleration ratio for residual-corrected low-precision quantization is only 3 or less(because there are other quantitative operations involved), often leading to less satisfactory acceleration results.

\section{Sparse quant}
Through the introduction in Section 2.3, there are two quantization methods to reduce errors, with the second method introducing a considerable amount of additional computation while reducing errors. Is there a way to minimize additional computation while reducing errors? In this section, we will introduce a method that introduces sparsity during the residual error correction process to reduce the overall computational load.

\subsection{Sparse residual error compensation algorithm}
\subsubsection{Methods to obtain sparsity}

From previous derivation , the residuals of $c_{ij}$ in the Equation(\ref{gemm_r}) result matrix $C$ can be expressed as (omits the multiplication terms of residuals with small values) :
\begin{equation}
\begin{split}
r{c^{ij}}=\sum_{k}((a_{ik}+Ra_{ik})((b_{kj}+Rb_{kj})) - a_{ik}b_{kj})\\
\approx \sum_{k}(ra_{ik}*b_{kj})+\sum_{k}(a_{ik}*rb_{kj})
\end{split}
\end{equation}

The relative error of $\delta{c_{ij}}$for each term of the result matrix $\delta c_{ij}$can be expressed as:

$$\delta{c_{ij}} \approx \sum_{k}ra_{ik}*b_{kj}/c_{ij}+\sum_{k}a_{ik}rb_{kj}/c_{ij}$$

Noticed that when performing low precision quantization:

\begin{eqnarray}    \label{eq}
    a^{Int}_{ij}=ToInt(\lambda *a^{Fp}_{ij})=\lfloor \lambda *a^{Fp}_{ij}\rfloor
\end{eqnarray}

$\lfloor.\rfloor$is the integer operation(We use downward rounding in our formulas to make it easier to understand). The multiplication operation is performed at the original floating point precision, so the error in this part is caused by the integer operation in Equation(\ref{countAint}).

\begin{eqnarray}    \label{countAint}
    Ra^{Int}_{ij}= \lfloor\lambda *a^{Fp}_{ij}\rfloor - \lambda *a^{Fp}_{ij} 
\end{eqnarray}

The transmission of this part of the error causes the final error of $a_{ij}$low-precision quantization, that is, the residual of $a_{ij}$, which can be expressed as:

\begin{eqnarray}    \label{eq}
    ra^{Fp}_{ij}&=&TypeCast(\lambda^{-1} *a^{Int}_{ij},Float) - a^{Fp} \nonumber \\
      ~&=& TypeCast((\lfloor\lambda *a^{Fp}_{ij}\rfloor-\lambda *a^{Fp}_{ij})
      *\lambda^{-1},Float) \nonumber \\
      ~&=& TypeCast((Ra^{Int}_{ij}*\lambda^{-1},Float)
\end{eqnarray}

Since $Ra^{Int}_{ij}\leq1$(Due to the rounding operation, the error will not exceed 1).
Thus 

\begin{equation}
\label{lambdaleq}
\begin{split}
ra^{Fp}_{ij}=TypeCast((Ra^{Int}_{ij}*\lambda^{-1},Float)\leq \lambda^{-1}
\end{split}
\end{equation}

Since Equation(\ref{lambdaleq}), $ra_{ik}$, $rb_{kj}$ is a number with a fixed range. Therefore, the impact on the entire $c_{ij}$ error mainly depends on the two terms $a_{ik} $, and $b_{kj}$:
$$\delta{c_{ij}}\leq
(\sum_{k}\lambda_a^{-1}b_{kj}/c_{ij}+\sum_{k}\lambda_b^{-1}a_{ik}/c_{ij})$$

That is, the relative error of the resulting matrix C is accumulated from the quantization exponents $\lambda_a$,$\lambda_b$and the original matrices A, B. That is to say, the numbers with larger values in A,B tend to have a larger impact on the residual matrix multiplication, while the relatively small numbers in the two matrices have a relatively small impact on the error of the resulting matrix C. Consider the absolute value of the numbers in the original matrix as the weights, that is, we only need to focus on the numbers with relatively large weights in the residual repair calculation to reduce the overall amount of calculation.

A simple way to do this is to set a threshold to remove the smaller values of the A and B matrices during the residual matrix multiplication, thus reducing the amount of computation. So how to measure the threshold?
Assume that we expect the relative error caused by matrix multiplication to the result $c_{ij}$to be at most 2M
$$\delta{c_{ij}}\leq (\sum_{k}\lambda_a^{-1}b_{kj}/c_{ij}+\sum_{k}\lambda_b^{-1}a_{ik}/c_{ij})\leq 2M$$

Each of the terms satisfied:

$$\lambda_a^{-1}\sum_{k}b_{kj}/c_{ij}\leq M,\lambda_b^{-1}\sum_{k}a_{ik}/c_{ij}\leq M$$

For each term in the summation:

\begin{equation}
    \begin{cases}
        a_{ik} \leq M*\lambda_b*(abs(c_{i1})/k) \\
        a_{ik} \leq M*\lambda_b*(abs(c_{i2})/k) \\ 
        ... \\
        a_{ik} \leq M*\lambda_b*(abs(c_{in})/k)
     \end{cases}
    \label{judgeofquant}
\end{equation}

Summarizing, A total of B.cols inequalities are satisfied:

\
$$\lambda_b^{-1}\sum_{k}a_{ik}/min(abs(c_{i*}))\leq M$$ 

Then, due to $c_{ij}$

$$a_{ik} \leq M*\lambda_b*(min(abs(c_{i*}))/k)$$

The elements of the A matrix satisfying this condition are put into the matrix $A^{'}$ and multiplied with the residual matrix $RB$ of $B$ to obtain the error compensation matrix, and the quantization result after the residual repair can be obtained

However, in practical experiments, it was observed that the scaling factor in the above formula was too large. By summing up the inequalities in Equation(\ref{judgeofquant}), the following scaling approach can result in a sparser matrix without significantly increasing the error.

    \begin{eqnarray}    \label{eq}
        a_{ik}\leq M*(avg(abs(c_{i1})))
        \label{avgQuan}
    \end{eqnarray}

where $c^{i*}$ can be estimated from the first quantized matrix multiplication result, without introducing additional computational overhead.

Through this formula, we can reduce the matrix A and B, and get the sparse matrix $A^{'},B^{'}$, and then reducing the calculation amount of low-precision residual compensation matrix multiplication.

Below is an example using the improved quantization method described above, adopting 8-bit quantization with $N=8$

\subsubsection{Example}
For input matrices A and B:
\[
A_{fp32} = \begin{bmatrix}
    0.2735 & -0.1588 & 0.1218 \\
    0.0953 & 1.5801 & -0.4861 \\
    -0.2394 & 0.1602 & 0.4294 \\
\end{bmatrix}
\]
\[
B_{fp32} = \begin{bmatrix}
    3.9284 & -0.0195 & -0.3836 \\
    -0.3288 & 2.2353 & -0.1895 \\
    -0.1376 & 0.0545 & -0.3641 \\
\end{bmatrix}
\]

The result matrix C for single-precision floating-point arithmetic is:
\[
C_{fp32} = \begin{bmatrix}
    1.1100 & -0.3538 & -0.1192 \\
    -0.0783 & 3.5038 & -0.1590 \\
    -1.0521 & 0.3861 & -0.0949 \\
\end{bmatrix}
\]

Compute the maximum values of matrices A and B, and determine the values of $\lambda_A,\lambda_B$ using Equation(\ref{getlambda})
with
\textbf{$max_A$} $= 1.58014$,
\textbf{$max_B$} $= 3.9284$,
then,
\textbf{$\lambda_A$} $= 40.5027$
\textbf{$\lambda_B$} $= 16.2916$

And compute the 8-bit quantized versions of matrices A and B:

\[
A_{int8} = \begin{bmatrix}
    11 & -6 & 4 \\
    3 & 64 & -19 \\
    -9 & 6 & 17 \\
\end{bmatrix} ,
B_{int8} = \begin{bmatrix}
    63 & 0 & -6 \\
    -5 & 36 & -3 \\
    -2 & 0 & -5 \\
\end{bmatrix}
\]
Compute matrix $C_{quant}$ by directly quantizing the result of low precision:
$C_{quant} = (A_{int8}*B_{int8})/(\lambda_A*\lambda_B)$

\[
C_{quant} = \begin{bmatrix}
    1.0836 & -0.3273 & -0.1031 \\
    -0.1409 & 3.4917 & -0.1743 \\
    -0.9563 & 0.3273 & -0.0743 \\
\end{bmatrix}
\]

Calculate the average values for each row and column of matrix $C_{quant}$

\[
C_{row} = \begin{bmatrix}
    0.5047 \\ 1.2690 \\ 0.4526 
\end{bmatrix} ,
C_{col}^T = \begin{bmatrix}
    0.7269 \\
    1.3821 \\
    0.1172 \\
\end{bmatrix}
\]

Dequantizing matrices $A,B$ and get $A^{'},B^{'}$
\[
A' = \begin{bmatrix}
    0.2716 & -0.1481 & 0.0988 \\
    0.0741 & 1.5801 & -0.4691 \\
    -0.2222 & 0.1481 & 0.4197 \\
\end{bmatrix}\\
\]
\[
B' = \begin{bmatrix}
    3.8670 & 0.0000 & -0.3683 \\
    -0.3069 & 2.2097 & -0.1841 \\
    -0.1228 & 0.0000 & -0.3069 \\
\end{bmatrix}
\]

Computing the residual matrix $RA,RB$ of matrix $A$ and $B$.
\[
RA = \begin{bmatrix}
    0.0020 & -0.0107 & 0.0231 \\
    0.0212 & 0.0000 & -0.0170 \\
    -0.0172 & 0.0120 & 0.0097 \\
\end{bmatrix}
\]
\[
Rb = \begin{bmatrix}
    0.0614 & -0.0195 & -0.0153 \\
    -0.0219 & 0.0256 & -0.0053 \\
    -0.0148 & 0.0545 & -0.0572 \\
\end{bmatrix}
\]

By Equation(\ref{avgQuan}) compute the reducing matrix $A^r,B^r$ from matrix $A$ and $B$.

\[
A^{r}  = \begin{bmatrix}
    0.2735 & 0.0000 & 0.0000 \\
    0.0000 & 1.5801 & 0.0000 \\
    0.0000 & 0.1602 & 0.4294 \\
\end{bmatrix}
\]

\[
B^{r} = \begin{bmatrix}
    3.9284 & 0.0000 & 0.0000 \\
    0.0000 & 2.2353 & 0.0000 \\
    0.0000 & 0.0000 & 0.0000 \\
\end{bmatrix}
\]
Quantize the matrices $A^{'} ,B^{'} $ in the same way to obtain their low-precision matrix:

\[
A_{int}^{'} = \begin{bmatrix}
    11 & 0 & 0 \\
    0 & 64 & 0 \\
    0 & 6 & 17 \\
\end{bmatrix},
B_{int}^{'} = \begin{bmatrix}
    63 & 0 & 0 \\
    0 & 36 & 0 \\
    0 & 0 & 0 \\
\end{bmatrix}
\]

Multiply the reduced matrices $A^{'} ,B^{'} $ with $RB,RA$ respectively. Add the results to $C_{quant}$ to obtain the matrix $C^{'}_{quant}$ after residual error correction:
$C^{'}_{quant} = C_q + A^{'} *RB+B^{'} *RA$

\[
C^{'}_{quant} = \begin{bmatrix}
    1.1077 & -0.3573 & -0.1070 \\
    -0.0878 & 3.5311 & -0.1819 \\
    -1.0356 & 0.3817 & -0.0987 \\
\end{bmatrix}
\]

Adopt the Frobenius norm of the matrix to estimate the absolute and relative errors of the result matrix. $E_r = \Vert X - X^* \Vert_{F} , E_\delta = \frac{E_r}{\Vert X \Vert_{F}} $
The quantization errors and direct quantization errors improved by the above algorithm are show in Table.\ref{tab:egresult}:

\begin{table}

  \caption{Errors of example}
    \begin{tabular}{|c|c|c|}
    \toprule
    Int8 quantization method & $E_r$    & $E_\delta$ \\
    \midrule
    Residual improvement & 0.144 & 0.011 \\
    \midrule
    Direct quantization & 0.799 & 0.035 \\
    \bottomrule
    \end{tabular}%
  \label{tab:egresult}%
\end{table}%

\subsubsection{Algorithm}
A simple example of a residual improvement algorithm is shown in the previous section. Besides,in practical execution, since the magnitude of sparsity s significantly impacts the execution speed of sparse matrix multiplication, in cases where the sparsity is not small enough, dense matrix multiplication exhibits faster execution speed.  As shown in Fig.\ref{FIG:diffper}(b), we can introduce a parameter $\eta$ (gemm has the same speed as spmm at $\eta$ sparsity) representing the acceleration ratio of sparse matrix multiplication to dense matrix multiplication at sparsity s. 

Only when $density<\eta$,  should the reduced matrix be used for sparse matrix multiplication. When $density>\eta$, it indicates that sparse matrix multiplication does not achieve significant acceleration within the given precision threshold. In such cases, the complete residual matrix multiplication is performed using the matrix before reduction. This ensures an optimal trade-off between maximizing quantization precision and minimizing computation time.

Algorithm.\ref{algorithmofigemm} summarizes the flow of the improved residual compensation algorithm. Firstly, the low precision matrix $D_{int}$ under N integer is calculated by direct quantization operation, and the original precision result $D_F$ is obtained by inverse quantization operation(line 1-3). Then the residual matrix $RA, RB$ of $A$ and $B$ is calculated by inverse quantization operation(line 4-6). And calculate D matrix for reduce matrix vector Drows, Dcols(line 7). The error threshold is used to obtain the possible sparse matrix $A^{'},B^{'}$ of matrix $A$ and $B$(line 8-10). And then select the appropriate matrix operation (GEMM,SPMM) under the sparsity of $A^{'},B^{'}$, and calculate the error compensation matrix D1,D2 of matrix $D$(line 11-15). Finally, $RD^1$ and $RD^2$ were used to compensate the error of the original quantization result matrix D and complete the remaining steps of GEMM(line 16-18).

\begin{algorithm}
    \label{algorithmofigemm}
  \SetAlgoLined
\KwData{$A,B,C$(Intputmatrix);$\alpha$,$\beta$(Scalar);\newline$\theta$(Threshold);$s$(Limit sparsity);\newline N(Quant bit),$T_F$(Origin precision )}
  
  \KwResult{matrix D}

    \tcc{Computes low-precision quantized matrix multiplication}
    $\{A_{int},B_{int}\}\leftarrow Quant(\{A,B\},N)$\;
    $D_{int} =$ GEMM$(A_{int},B_{int})$\; \tcc{compute $D = A_{int}\cdot B_{int}$}
    $D_{F} \leftarrow Dequant(D_{int},T_F)$\;

    \tcc{Calculate the residual matrix}
    $\{A_{F},B_{F}\} \leftarrow Dequant(\{A_{int},B_{int}\},T_F)$\;
    $\{RA,RB\} \leftarrow \{A,B\} - \{A_{F},B_{F}\}$\;
    $\{RA_{int},RB_{int}\} \leftarrow Quant(\{RA,RB\},N)$\;

    \tcc{Reduce sparse matrix}
    $\{D^{row},D^{col}\} = GetAvgVec(D_{F},\{row,col\})$ \tcp*[r]{compute avg vector of cols and rows from $D_{F}$}
    $\{A^{'},B^{'}\} \leftarrow Reduce(\{A,B\},\{D^{row},D^{col}\},\theta)$ \tcp*[r]{reduce sparse matrix from A and B}
    $\{sp_A,sp_B\} = getSparsity(\{A^{'},B^{'}\})$ \;
    $\{A^{'}_{int},B^{'}_{int}\}\leftarrow Quant(\{A^{'},B^{'}\},N)$\;

    \tcc{Computes residual matrix multiplication}
    \eIf{$\{sp_A,sp_B\}<s$}{
        $\{DR^{1}_{int},DR^{2}_{int}\}$=SPMM($\{A^{'},B^{'}\},\{RB,RA\}$)\;
    } {
        $\{DR^{1}_{int},DR^{2}_{int}\}$ $=$ GEMM($\{A,B\},\{RB,RA\}$)\;
    }

    \tcc{Error compensation}
    $\{RD^{1}_{F},RD^{2}_{F}\} \leftarrow Dequant(\{RD^{1}_{int},RD^{2}_{int}\},T_F)$\;
    $D_{F} = D_{F} + RD^{1}_{F} + RD^{2}_{F}$\;

    \Return $D_{F} = \alpha*D_{F} + \beta C$\;
    
  \caption{Algorithms of compute GEMM $D=\alpha A\cdot B+\beta C$ in Int-quantization }
\end{algorithm}

\subsection{Vector-wise Quantization}

Through the above algorithm analysis, the error caused by the proposed method to the final structure consists of three parts $Er_{sum} = Er+Es + Eq$: 
Where $Er_{mul} = A^r*B^r$ is the omitted residual multiplication term; Multiplication of sparse omitted item ${Es} =\sum_{k}ra_{ik}*b_{kj}+\sum_{k}a_{ik}*rb_{kj}$ (for $a_{ik},b_{kj}$ satisfy reduce conditions); Error $E_q$caused by quantizing the residual term.

For these three error terms, improving the quantization-induced error in the third term by applying iterative residual compensation can enhance the accuracy. However, this approach introduces greater computational complexity, and the resulting precision improvement may not justify the costs \cite{ref31}. In contrast, the first two error terms, if minimized due to the quantization-induced error being smaller, lead to a smaller value for $A^r * B^r$. Consequently, the error introduced to the result is also smaller.

\begin{figure*}
    \centering
    \includegraphics[width=1\linewidth]{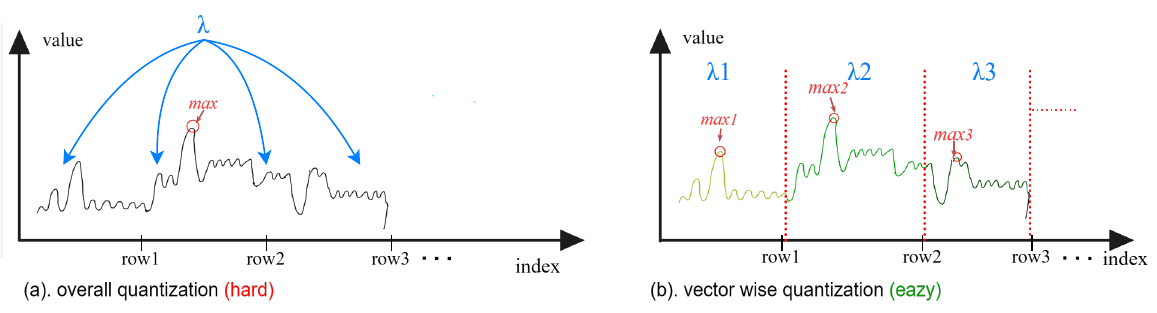}
    \caption{Comparison of the effect of different quantization methods}
    \label{Fig:Comparison_methods}
\end{figure*}

As shown in Figure.\ref{Fig:Comparison_methods}. When the scale of matrix multiplication is larger, the values within the matrix exhibit a broader range of variations, making quantization susceptible to the influence of extreme values and causing significant errors. In contrast, the range of variations for each row and column within the matrix is smaller. Therefore, a natural approach to reduce errors is to introduce a method based on row-wise and column-wise quantization. By selecting the maximum value relative to each row and column, the error is minimized. This method balances the difficulty of quantization for both matrices.

The quantization process using this method follows these steps: first, instead of globally selecting the absolute maximum value for the entire matrix, the absolute maximum value is chosen for each row of matrix A and each column of matrix B. Then, the corresponding $\lambda$ values for each row and column of A and B are calculated using the formula. Subsequently, the quantization process is applied. Finally, during the dequantization of the computed low-precision matrix $C_{fp} = \frac{C_{int}}{\lambda_{a_i} \times \lambda_{b_j}}$.

\subsection{Estimation of Algorithm Execution Flow Time Complexity}
In this section, we delve into the time complexity estimation of the proposed algorithm's execution flow. By analyzing the key steps and operations involved, we aim to provide insights into the computational efficiency of our approach.

The execution flow of the improved residual repair algorithm is as follows, which can be represented in Figure 6:

In order to show the superiority of the algorithm, we take the size of the data into consideration in the calculation of time complexity, assuming that the original precision is P bit, and the low precision quantization is Q bit. The input and output matrices are both of size $N*N$.
The improved residual repair matrix multiplication has three parts: normal low-precision quantized matrix multiplication, sparse reduction, and residual matrix multiplication
\begin{enumerate}
    \item Normal low-precision quantized matrix multiplication: The initial stage involves quantizing and dequantizing matrices A and B and GEMM computing.
    \begin{itemize}
        \item Firstly, $A$ and $B$ matrices are vector-wise quantized, and the quantized matrix $A_{int}$,$B_{int}$ - $O(2*(PN)^2)$
        \item And then through $A_{int}$, $A_{int}$ matrix multiplication to get the quantized result matrix - $O((QN)^3)$
        \item Deuantization of the vector-wise quantization result matrix Cint is performed to obtain the calculation result - $O((PN)^2)$
    \end{itemize}
    \item Sparse reduction: In this step, the residual matrix quantified in the previous step is calculated and the sparse residual matrix is obtained by removing the numbers with less weight using the reduction formula.
    \begin{itemize}
        \item The residual matrix $AR$, $BR$ is obtained by inverse quantization of the quantization matrix  $A_{int}$,$B_{int}$ - $O(2*(PN)^2)$
        \item The residual matrix is quantized to get ARint, BRint - $O(2*(PN)^2)$
        \item A and B matrices are reduced by using the above results, and A possible sparse matrix $A'$, $B'$ is obtained - $O(2*(PN)^2)$
    \end{itemize}    
    \item Residual matrix multiplication: In this step, the reduced sparse matrix multiplication is performed and the result is inversely quantized and added to the original quantization result.
    \begin{itemize}
        \item A', B' is multiplied by the residual matrix, and calculated $CR_{int}$ .The time complexity for sparse matrix multiplication can be approximated as $O(2S*(QN)^3)$, where S represents the spmm acceleration ratio.
        \item Dequantization of  $CR_{int}$ is performed to obtain the residual error compensation matrix CR, and the time complexity is$O((PN)^2)$
        \item Use $C = C + CR + CR2$for error compensation - $O((PN)^2)$
    \end{itemize}    
\end{enumerate}

In addition, the total time complexity of the improved algorithm is: $11*O((PN)^2)+ (1+2S)*O((QN)^3)$

The total time complexity of the original low precision quantization algorithm is: $3*O((32n)^2)+ O((16n)^3)$

If fullsize's residual matrix repair method is used, the Sparse reduction in step 2 is not required, but the time complexity of the residual matrix multiplication is $O(2*(QN)^3)$. The total time complexity of the algorithm is: $9*O((PN)^2)+ 3*O((QN)^3)$

To sum up, compared with fullsize's residual matrix repair method, there is a significant reduction about three times in the calculation of the cubic time complexity if the reduce matrix sparsity is high enough.

\begin{figure}
	\centering
\includegraphics[width=1\linewidth]{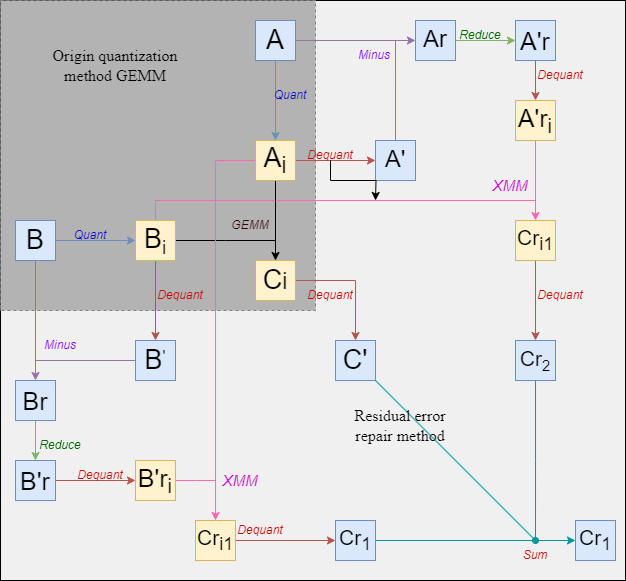}
	\caption{Workflow of Algorithm}
	\label{FIG:6}
\end{figure}

\section{Evaluation}

The experimental setup involved the following computing devices:

\textbf{CPU:}
\begin{itemize}
    \item Model: Intel(R) Xeon(R) Platinum 8163 CPU
    \item Compilation:  \texttt{-O3 -mavx,-mfma,-march=native}
    \item Operating System: Ubuntu 11.4.0-1ubuntu1~22.04
\end{itemize}

\textbf{GPU:}
\begin{itemize}
    \item Model: NVIDIA A100-PCIE-40GB
    \item Compilation: \texttt{nvcc} with options \texttt{-O3, -gencode \\arch=compute\_80, code=sm\_80}
\end{itemize}

The algorithm implementation, as described in Section 3.3 , consists of three main components:

\begin{itemize}

    \item \textbf{Dense Low-Precision Matrix Multiplication}: \texttt{CUTLASS (version 3.3)} on GPU and \texttt{MKL} on CPU was utilized to construct the matrix multiplication kernel.
    \item \textbf{Sparse Low-Precision Multiplication}: The \texttt{cusparse SpMM} function on GPU and \texttt{MKL} on CPU was employed for sparse matrix multiplication. \textit{(Since \texttt{MKL} does not support sparse matrix operations that precision below Int32, and we did not find a  high-performance low-precision sparse math library on the CPU. Instead, Int32 SPMM is used for this part, and speed estimates are estimated by Int8 GEMM in \texttt{MKL}.)}
    \item \textbf{Matrix Quantization and Reduction Operations}: Developed using the \texttt{CUDA toolkit (version 12.3)} on GPU and \texttt{OpenMP} on CPU for matrix quantization and reduce operations.
    
\end{itemize}

\subsection{GEMM Effectiveness}

Through three distinct evaluations, we assess the effectiveness of the improved residual error compensation algorithm in comparison to matrix multiplication. These evaluations encompass precision testing under varying thresholds and matrix sizes, runtime testing at different sparsity levels, and an analysis of the contribution of quantization strategies throughout the algorithm.

\subsubsection{Precision Testing}
Precision testtest is carried out using single precision SGEMM as reference. The error metric $E_r$, introduced earlier, was utilized for assessment. The error results obtained using this method for matrix multiplication with different distributions are presented in the following Fig.\ref{Fig:rxigemm-precision}

Figure observations indicate that the algorithm employing a hybrid strategy can achieve considerable computational precision below a certain threshold level. Compared to the full residual quantization approach, our hybrid strategy algorithm typically attains higher accuracy at lower threshold values. Moreover, our hybrid quantization approach yields a substantial precision increase when benchmarked against the traditional direct quantization method, achieving at least an 80\% reduction in relative error under int8 quantization with a chi-squared distribution.

Moreover, due to the incorporation of the vector-wise method, our algorithm still demonstrates significant advantages under int4. As shown in Fig.\ref{Fig:rxigemm-precision}(a,c), Compared with the fully residual quantization method, the improved algorithm can obtain higher accuracy than the complete residual quantization method at a higher threshold (th=0.5).

It is noteworthy that the results of the hybrid strategy at lower precisions (int4) do not vary with the threshold when assessing exponential and chi-squared distributions because the preponderance of data is clustered near smaller magnitudes with rarified instances of markedly larger values. Hence, increasing the threshold invariably preserves these exceptional values within the sparse matrix reduction step, maintaining consistent relative errors. Simultaneously, precision improves more significantly as the threshold decreases. The hybrid strategy, even at reduced precision levels, can surpass the computational accuracy of the full residual method by nearly 10\% under certain distributions, achieving better computational precision at elevated threshold levels (e.g. Th=0.2).

Under int8 quantization, the algorithm can still achieve good calculation accuracy. With the reduction of the threshold, the calculation accuracy will also be significantly improved, and eventually achieve full residual quantization, and even higher accuracy. Especially in the Chi-square distribution, the improved algorithm can achieve at least 5 times the accuracy improvement of the original quantization method, and can achieve similar accuracy to the full residual method under a fairly high threshold (th=1)

\begin{figure}
    \centering
    \includegraphics[width=1\linewidth]{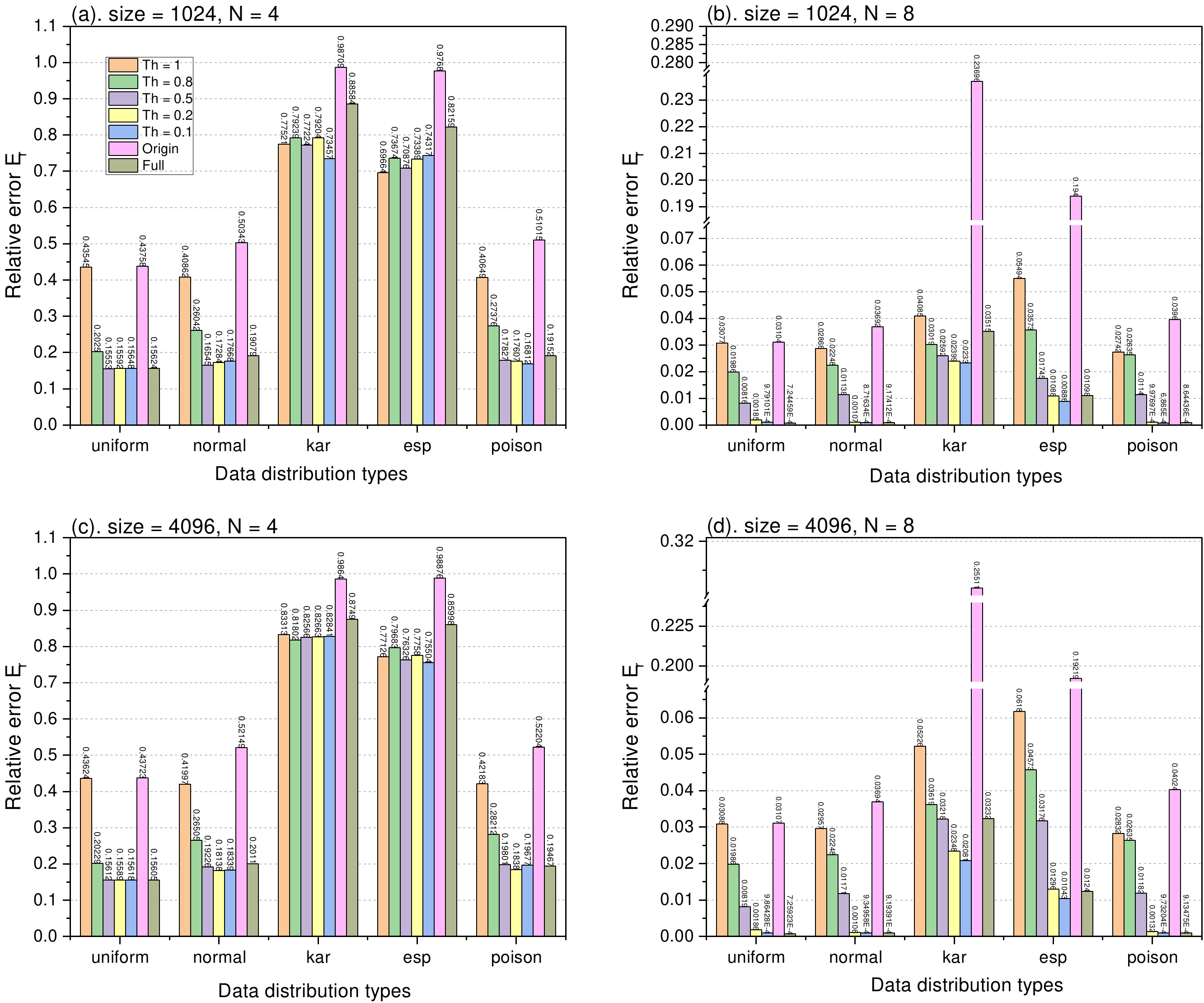}

    \caption{Precision test of Algorithms(int4/int8,size=1024), where \textbf{'Full'} represents  residual quantization matrix multiplication, and \textbf{'origin'} represents the original method of direct quantization and \textbf{'Th'} denotes the hyperparameter of the sparse quantization algorithm mentioned earlier—threshold. Where \textbf{'uniform'} is a standard random distribution of 0 to 1.   \textbf{'Normal'} is a normal distribution with a mean of 10 and a variance of 3. \textbf{'Esp'} is the exponential distribution of $\lambda=4$. \textbf{'Poison'} is the Poisson distribution of $\lambda=10$. \textbf{'Kar'} is a chi-square distribution of n=1.}
    \label{Fig:rxigemm-precision}
\end{figure}

\subsubsection{Quantization Strategy Testing}
\begin{figure}
    \centering
    \includegraphics[width=1\linewidth]{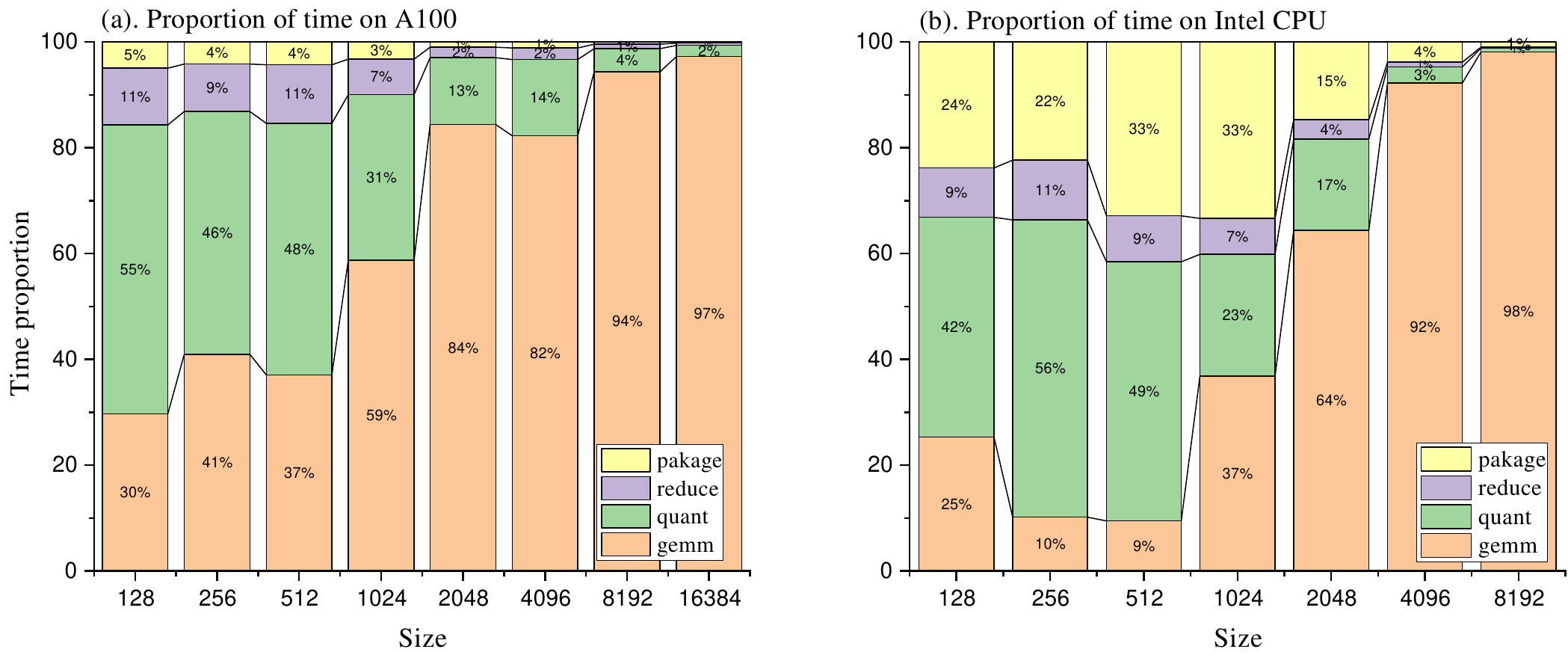}
    \caption{Time proportion of different parts of the algorithm}
    \label{Quantization Strategy Testing}
\end{figure}

Figure.\ref{Quantization Strategy Testing} illustrates the time distribution of various quantization steps for different matrix sizes on diverse computing platforms. The algorithm prioritizes operations with lower relative execution times (spmm/gemm) during residual error compensation. The steps include quantization/unquantization (quant), low-precision matrix multiplication (xxmm), sparse matrix reduction (reduce), and the time consumed by residual matrix error compensation and other components (package). The analysis assumes worst-case scenario conditions where the algorithm operates under the premise that the total time for quantization equals that of gemm. Fig.\ref{Quantization Strategy Testing} indicates that, with increasing matrix size, gemm/spmm occupy a larger share of the total time (approximately 90\%). The additional computation introduced by quantization strategies remains acceptable.

\subsubsection{Speed Testing}
In light of the considerable efficiency discrepancy between Nvidia-cuSparse and cublas on Nvidia A100 for Sparse Matrix-Matrix Multiplication (SPMM) and General Matrix-Matrix Multiplication (GEMM), it is observed that the performance of sparse matrix multiplication only surpasses the predefined threshold $\eta$ (approximately 0.001), indicating that the enhanced algorithm resorts to SPMM exclusively when the matrix density falls below 0.001. Consequently, the velocity of the advanced quantized residual algorithm is assessed on the Nvidia A100, using int8 low-precision multiplication, across the density spectrum of 0.001 to 0.00005. But on the CPU platform, sparse matrix can achieve a high acceleration, $\eta$(approximately 0.3), so the matrix density used in the experiment ranges from 0.3 to 0.01.The execution employs matrices constituted of uniformly distributed random values ranging from 0 to 1. Table \ref{speed test gpu} and Table \ref{speed test cpu} presents the architecture of the improved algorithm's integer-floating point computational speed in contrast with two benchmark methods on GPU and CPU.

The tabulated data reveals an appreciable acceleration effect by the refined residual patching technique as matrix density diminishes. On the two computing devices,  the enhanced approach attains an acceleration rate approximately 2.5 times greater than the complete residual as sparsity increases. This underscores the efficacy of the refined algorithm.

\begin{table*}[htbp]
  \centering
  \caption{\centering{TOPS of algorithm at different densities on Nvidia-A100}}
    \captionsetup{justification=centering, labelsep=newline} 
    \begin{tabular}{|c|c|c|c|c|c|r|c|}
    \toprule
    \multirow{2}[4]{*}{Size} & \multicolumn{5}{c|}{TOPS to density on GPU} & \multicolumn{2}{c|}{\multirow{2}[4]{*}{Baseline methods}} \\
\cmidrule{2-6}          & 0.001 & 0.0005 & 0.00025 & 0.0001 & 0.00005 & \multicolumn{2}{c|}{} \\
    \midrule
    \hline
    \multirow{2}[4]{*}{1024} & \multirow{2}[4]{*}{3.21} & \multirow{2}[4]{*}{3.61} & \multirow{2}[4]{*}{4.22} & \multirow{2}[4]{*}{5.31} & \multirow{2}[4]{*}{5.62} & \multicolumn{1}{c|}{3.35} & Full \\
\cmidrule{7-8}          &       &       &       &       &       & 9.06 & Origin \\
    \midrule
    \multirow{2}[4]{*}{4096} & \multirow{2}[4]{*}{25.19} & \multirow{2}[4]{*}{32.36} & \multirow{2}[4]{*}{44.13} & \multirow{2}[4]{*}{53.13} & \multirow{2}[4]{*}{68.12} & \multicolumn{1}{c|}{27.96} & Full \\
\cmidrule{7-8}          &       &       &       &       &       & 83.97 & Origin \\
    \midrule
    \multirow{2}[4]{*}{8192} & \multirow{2}[4]{*}{40.33} & \multirow{2}[4]{*}{49.9} & \multirow{2}[4]{*}{67.72} & \multirow{2}[4]{*}{74.48} & \multirow{2}[4]{*}{97.35} & \multicolumn{1}{c|}{44.05} & Full \\
\cmidrule{7-8}          &       &       &       &       &       & 125.06 & Origin \\
    \midrule
    \multirow{2}[4]{*}{16384} & \multirow{2}[4]{*}{49.82} & \multirow{2}[4]{*}{60.12} & \multirow{2}[4]{*}{77.73} & \multirow{2}[4]{*}{100.41} & \multirow{2}[4]{*}{128.96} & \multicolumn{1}{c|}{50.32} & Full \\
\cmidrule{7-8}          &       &       &       &       &       & 145.20 & Origin \\
    \bottomrule
    \end{tabular}%

  \label{speed test gpu}%
\end{table*}%

\begin{table*}[htbp]
  \centering
  
  \caption{\centering{GOPS of algorithm at different densities on Intel(R) Xeon(R) Platinum 8163 CPU}}
  \captionsetup{justification=centering, labelsep=newline} 
    \begin{tabular}{|c|c|c|c|c|c|c|c|}
    \toprule
    \multirow{2}[4]{*}{Size} & \multicolumn{5}{c|}{GOPS to density on CPU} & \multicolumn{2}{c|}{\multirow{2}[4]{*}{Baseline methods}} \\
\cmidrule{2-6}          & 0.30  & 0.20  & 0.10  & 0.05  & 0.01  & \multicolumn{2}{c|}{} \\
    \midrule
    \hline
    \multirow{2}[4]{*}{512 } & \multirow{2}[4]{*}{7.84 } & \multirow{2}[4]{*}{7.87 } & \multirow{2}[4]{*}{7.89 } & \multirow{2}[4]{*}{7.92 } & \multirow{2}[4]{*}{7.92 } & 9.19  & Full \\
\cmidrule{7-8}          &       &       &       &       &       & 23.39  & Origin \\
    \midrule
    \multirow{2}[4]{*}{1024 } & \multirow{2}[4]{*}{67.28 } & \multirow{2}[4]{*}{70.05 } & \multirow{2}[4]{*}{71.32 } & \multirow{2}[4]{*}{71.91 } & \multirow{2}[4]{*}{72.24 } & 73.56  & Full \\
\cmidrule{7-8}          &       &       &       &       &       & 175.76  & Origin \\
    \midrule
    \multirow{2}[4]{*}{4096 } & \multirow{2}[4]{*}{290.44 } & \multirow{2}[4]{*}{354.21 } & \multirow{2}[4]{*}{443.84 } & \multirow{2}[4]{*}{513.31 } & \multirow{2}[4]{*}{582.08 } & 277.40  & Full \\
\cmidrule{7-8}          &       &       &       &       &       & 799.11  & Origin \\
    \midrule
    \multirow{2}[4]{*}{8192 } & \multirow{2}[4]{*}{335.83 } & \multirow{2}[4]{*}{469.07 } & \multirow{2}[4]{*}{572.40 } & \multirow{2}[4]{*}{717.81 } & \multirow{2}[4]{*}{926.42 } & 379.75  & Full \\
\cmidrule{7-8}          &       &       &       &       &       & 1122.17  & Origin \\
    \bottomrule
    \end{tabular}%
  \label{speed test cpu}%
\end{table*}%

\begin{figure}
    \centering
    \includegraphics[width=1\linewidth]{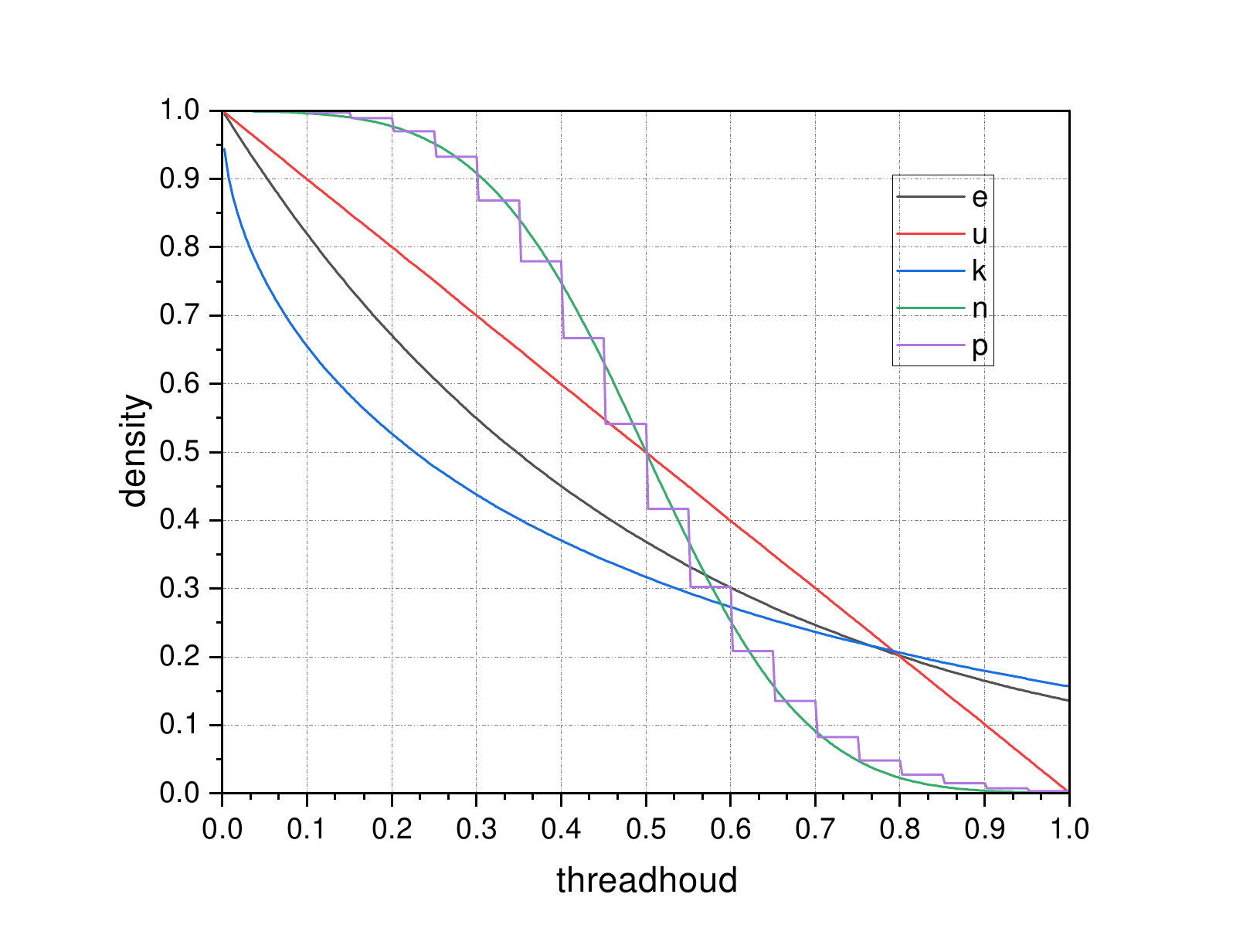}

    \caption{The reduced matrix density of matrices generated by different random methods.}
    \label{Fig:density}
\end{figure}

\subsubsection{Comparison and Conclusion}
Through the experiments of the above three parts, we understand the improved low precision quantization algorithm respectively from the three aspects of speed, algorithm structure and accuracy. In order to further understand the effectiveness and limitations of the algorithm, we present the Fig.\ref{Fig:density} of test matrix density with different distributions used in the experiment along with the threshold value.

Since the sparse matrix execution efficiency of GPU is generally low, the improved algorithm is difficult to show enough reduce sparsity in the naturally randomly generated matrix to obtain good acceleration effect. We will analyze the performance test structure on the CPU platform.

In Fig.\ref{FIG:diffper}, it shows that $\eta$ is about 0.3 on the CPU platform. Combined with Fig.\ref{Fig:density}. Then it can derive the thresholds at which the acceleration effect can occur under different distributions: In Esp(4),Normal(10,3) and Poison(10), it is about 0.6; In Kar(1) is about 0.5. In Uniform(0,1) is about 0.7. In other words, under these distributions, only when the threshold value is greater than these values can sparse computation be accelerated.

As shown in Fig.\ref{Fig:rxigemm-precision}, in the data of int4, combined with the analysis of the previous paragraph, the uniform distribution can be obtained at TH=0.8, and the calculation acceleration of sparse matrix can be obtained at 1. When TH=0.8, the matrix density is about 0.2. Combined with Table.\ref{speed test cpu},  the calculation acceleration of about 1.2 times can be obtained at this density compared with the complete residual quantization method.  Using this method, the parts of the algorithm that generate sparsity acceleration are analyzed, and the complete residual quantization method is compared. The results obtained are shown in table.\ref{accuracy and speed}, where the green part indicates that both accuracy and speed have been improved, the black part indicates that the accuracy has only been accelerated without accuracy improvement, and the blue part indicates that the accuracy has hardly changed, but the performance has been improved. The red part indicates that the sparsity of the threshold value is not enough to call SPMM for acceleration.

The algorithm can achieve significant results in some distributions, and the algorithm can achieve performance improvements of up to 15\% and 45\% in esp distributions under int4. In the case of int8, the algorithm can achieve a similar accuracy to the full residual with a threshold value of 0.8 under the kar square distribution, but the speed is increased by 28\%.

In other cases, although the accuracy of the algorithm is not as good as that of the complete residual, compared with the direct quantization method, the efficiency of the improved algorithm can be close to that of the direct quantization method when the reduce matrix is sparse enough. In conclusion, the improved algorithm is an intermediate method that can measure efficiency by threshold value.

\begin{table*}[htbp]
  \centering
  \caption{Improved accuracy and speedup under different thresholds on CPU}
    \begin{tabular}{|c|c|c|c|c|c|cc|}
    \toprule
    \multicolumn{8}{|c|}{ percision enhancement and Spedup at threadhoud} \\
    \midrule
    \multirow{2}[4]{*}{Quant type} & random type & \multicolumn{2}{c|}{Th = 1} & \multicolumn{2}{c|}{Th = 0.8} & \multicolumn{2}{c|}{Th = 0.5} \\
\cmidrule{2-8}          &       & acc  up & speedup & acc  up & speedup & \multicolumn{1}{c|}{acc  up} & speedup \\
    \midrule
    \multirow{5}[10]{*}{xigemm-int8} & uniform & -41.63  & 2.10  & -26.41  & 1.28  & \multicolumn{2}{c|}{\textcolor[rgb]{ 1,  0,  0}{density insufficient}} \\
\cmidrule{2-8}          & normal & -31.23  & 2.10  & -23.50  & 1.85  & \multicolumn{2}{c|}{\textcolor[rgb]{ 1,  0,  0}{density insufficient}} \\
\cmidrule{2-8}          & kar   & -0.49  & 1.38  & \textcolor[rgb]{ 0,  .69,  .941}{-0.03 } & \textcolor[rgb]{ 0,  .69,  .941}{1.28 } & \multicolumn{1}{c|}{\textcolor[rgb]{ 0,  .69,  .314}{0.08 }} & \textcolor[rgb]{ 0,  .69,  .314}{1.05 } \\
\cmidrule{2-8}          & esp   & -4.63  & 1.46  & -2.25  & 1.24  & \multicolumn{2}{c|}{\textcolor[rgb]{ 1,  0,  0}{density insufficient}} \\
\cmidrule{2-8}          & poison & -31.76  & 2.10  & -29.48  & 1.95  & \multicolumn{2}{c|}{\textcolor[rgb]{ 1,  0,  0}{density insufficient}} \\
    \midrule
    \multirow{5}[10]{*}{xigemm-int4} & uniform & -1.79  & 2.10  & -0.30  & 1.28  & \multicolumn{2}{c|}{\textcolor[rgb]{ 1,  0,  0}{density insufficient}} \\
\cmidrule{2-8}          & normal & -1.14  & 2.10  & -0.36  & 1.85  & \multicolumn{2}{c|}{\textcolor[rgb]{ 1,  0,  0}{density insufficient}} \\
\cmidrule{2-8}          & kar   & \textcolor[rgb]{ 0,  .69,  .314}{0.12 } & \textcolor[rgb]{ 0,  .69,  .314}{1.38 } & \textcolor[rgb]{ 0,  .69,  .314}{0.11 } & \textcolor[rgb]{ 0,  .69,  .314}{1.28 } & \multicolumn{1}{c|}{\textcolor[rgb]{ 0,  .69,  .314}{0.13 }} & \textcolor[rgb]{ 0,  .69,  .314}{1.05 } \\
\cmidrule{2-8}          & esp   & \textcolor[rgb]{ 0,  .69,  .314}{0.15 } & \textcolor[rgb]{ 0,  .69,  .314}{1.46 } & \textcolor[rgb]{ 0,  .69,  .314}{0.10 } & \textcolor[rgb]{ 0,  .69,  .314}{1.24 } & \multicolumn{2}{c|}{\textcolor[rgb]{ 1,  0,  0}{density insufficient}} \\
\cmidrule{2-8}          & poison & -1.12  & 2.10  & -0.43  & 1.95  & \multicolumn{2}{c|}{\textcolor[rgb]{ 1,  0,  0}{density insufficient}} \\
    \bottomrule
    \end{tabular}%
  \label{accuracy and speed}%
\end{table*}%

The comprehensive testing above reveals that the efficiency of low-precision matrix operations, coupled with precision control through error thresholds and the judicious selection of computational methods, ensures effective low-precision acceleration while maintaining required precision levels. Moreover, the adopted quantization optimization exhibits a decreasing contribution to the overall computation time as matrix size increases, affirming the effectiveness of the quantization strategy and the superiority of the overall algorithm.

\subsection{QR decomposition}
In the QR decomposition routine employing Householder decomposition, the application of the orthogonal matrix to the upper triangular result matrix can be substituted with the improved sparse matrix quantized multiplication.In the presented analysis, the relative error of standard multiplication is defined as the relative error between the matrix A' obtained by direct single-precision floating-point multiplication (A' = Q * R) and the original matrix A. The subsequent three cases involve different quantization methods, each compared in terms of relative error against the standard A matrix.The results are shown in Table \ref{method in qr}.

It is evident that employing sparse matrix precision compensation yields a noticeable improvement in accuracy compared to direct quantization. The matrices exhibit a considerable level of sparsity, which is notably influenced by quantization thresholds, bit precision, and the original matrix shape. As the matrix size increases, the step involving the multiplication of orthogonal matrices in the Householder transformation results in a markedly sparse matrix (with sparsity around 1e-3 under conditions of M=N=K, threshold of 1e-4, and 16-bit quantization).
When utilizing the Frobenius norm of the matrix as the relative error metric, defined as $E_r = \Vert X - X^* \Vert_{F} / \Vert X \Vert_{F} $
the errors are as follows for matrices generated from a normal distribution:

In the provided context, FP32 represents the result of single-precision floating-point arithmetic. Quan-ori, Quan-res, and Quan-spres respectively denote the original single-precision matrix multiplication, the complete residual error compensation algorithm, and the improved sparse residual error compensation algorithm.

The table clearly indicates that the improved residual error compensation algorithm achieves a comparable error magnitude to dense quantization while simultaneously reducing computational complexity.

\begin{table}[htbp]
  \centering
  \caption{Different mixing precision methods in QR decomposition}
    \begin{tabular}{|c|c|c|c|}
    \toprule
    \multirow{2}[4]{*}{matrix type} & \multirow{2}[4]{*}{method} & int4  & int8 \\
\cmidrule{3-4}          &       & ER    & ER \\
    \midrule
    \multirow{3}[6]{*}{kar} & origin & 0.945681 & 0.650822 \\
\cmidrule{2-4}          & full  & 0.933172 & 0.121827 \\
\cmidrule{2-4}          & xigemm & 0.930159 & 0.150565 \\
    \midrule
    \multirow{3}[6]{*}{uniform} & origin & 0.979227 & 0.819087 \\
\cmidrule{2-4}          & full  & 0.97261 & 0.139795 \\
\cmidrule{2-4}          & xigemm & 0.96432 & 0.141484 \\
    \midrule
    \multirow{3}[6]{*}{normal} & origin & 0.995767 & 0.835512 \\
\cmidrule{2-4}          & full  & 0.99078 & 0.124808 \\
\cmidrule{2-4}          & xigemm & 0.993823 & 0.146129 \\
    \midrule
    \multirow{3}[6]{*}{esp} & origin & 0.938367 & 0.669409 \\
\cmidrule{2-4}          & full  & 0.932157 & 0.121116 \\
\cmidrule{2-4}          & xigemm & 0.926652 & 0.149053 \\
    \bottomrule
    \end{tabular}%
  \label{method in qr}%
\end{table}%

\subsection{Convolutional neural network}

LeNet-5 is a convolutional neural network consisting of 3 convolutional layers, 2 pooling layers, and 3 fully connected layers. In this experiment, a modified version of LeNet-5 with reduced layers is employed, containing one convolutional layer and two fully connected layers. Various quantized matrix multiplication techniques proposed in this study are applied specifically in the fully connected layers. The convolutional operation is processed using im2col, while the fully connected layers are transformed into matrix multiplication through concatenation.

The experimental results are presented below, where 'n' represents the quantization bit width, and 'th' denotes the hyperparameter of the sparse quantization algorithm mentioned earlier—threshold. The dataset utilized is FashionMNIST, comprising 10,000 images, with a batch size of 10,000.

The precision and speedup ratios obtained with different quantization thresholds are detailed in the table below, which is evident that the improved algorithm achieves a substantial speedup.

\section{Discussion}
\subsection{Performance evaluation}
The efficiency of the proposed algorithm depends on two main aspects:

\textbf{Efficiency of Low-Precision Operations:}
Performing tensor computations with low-precision integer types not only achieves computational acceleration with acceptable loss but also significantly saves energy.

\textbf{Calculation Efficiency of Sparse Matrix Multiplication compared to Dense Matrix Multiplication} :
The efficiency of sparse matrix operations varies depending on the form of sparsity. In this work, the sparse matrix is represented in Compressed Sparse Row (CSR) format. It's worth noting that different sparse matrix formats, such as Block Sparse Row (BSR) format for block sparse matrices, might exhibit better performance, especially when utilizing tensor cores. This aspect, however, is not explored further in this paper. Sparse matrices, especially when dealing with matrices containing a low number of non-zero elements, demonstrate significant performance improvements in matrix computations. Additionally, the sparsity of the generated sparse matrix is determined by the error threshold and the structure of the input matrix. When the input matrix contains larger values, the resulting residual correction matrix tends to be sparser.

Besides, Nvidia supports dense matrix multiplication with sparsity starting with the Ampare architecture (using mma as the multiplication core instruction). In the future, we can explore the possibility of using sparse dense matrix multiplication instead of Spmm when sparsity is large to obtain a good acceleration ratio.

In summary, compared to fully residual correction methods, the proposed approach yields performance improvements when the sparsity is sufficiently high.

\begin{table*}[h]
  \centering
  \caption{The acceleration effect of the algorithm in the neural network}
    \begin{tabular}{|c|c|c|c|c|c|c|}
    \toprule
    \multirow{8}[16]{*}{Lenet-5} & \multirow{3}[6]{*}{int4} & th    & 1     & 0.1   & 0.01  & 0 \\
\cmidrule{3-7}          &       & acc   & 0.4294 & 0.5793 & 0.671 & 0.7023 \\
\cmidrule{3-7}          &       & speedup & 10.13x & 6.62x & 5.14x & 3.55x \\
\cmidrule{2-7}          & \multicolumn{6}{c|}{} \\
\cmidrule{2-7}          & \multirow{3}[6]{*}{int8} & th    & 1     & 0.1   & 0.01  & 0 \\
\cmidrule{3-7}          &       & acc   & 0.8403 & 0.8486 & 0.8512 & 0.8518 \\
\cmidrule{3-7}          &       & speedup & 2.32x & 1.96x & 1.42x & 1.42x \\
\cmidrule{2-7}          & float & acc   & \multicolumn{4}{c|}{0.8518} \\
    \bottomrule
    \end{tabular}%
  \label{tab:addlabel}%
\end{table*}%

\subsection{Hyperparameter fine-tuning}

The algorithm proposed in this paper involves two hyperparameters: the error threshold $M$ and the ratio of the execution speeds of GEMM (General Matrix Multiply) to SPMM (Sparse Matrix Multiply). The latter is a constant value specific to each computing device and can be determined before the algorithm starts.

The error threshold represents an expected value for the residual correction error in most cases, serving as an upper limit for residual correction errors. However, due to the scaling of the mean, the error may exceed the threshold in the presence of some singular values (extremely large or small values). Additionally, considering the limitations of low precision, the minimum error achievable with the residual compensation method will not surpass the error of the fully residual compensation method. The practical minimum error for the error threshold corresponds to the error of the fully residual compensation method.

To strike a suitable balance between performance and accuracy for application implementation, a functionality for automatically adjusting the error threshold will be developed. This will allow fine-tuning the error threshold to achieve the desired trade-off between precision and performance.

\subsection{Possible application}

Low-precision quantization has gained widespread application in the field of machine learning. The proposed improved residual compensation in this paper demonstrates significant accuracy improvements with minimal computational cost in certain scenarios. In numerical computing, where algorithms are sensitive to errors, the residual compensation method presented in this paper achieves substantial error reduction in low precision. However, the error magnitudes still fall short of meeting mathematical requirements, posing a significant challenge for low-precision applications in traditional High-Performance Computing (HPC).

The GMRES-IR method for solving LU equations using iterative techniques introduces low-precision gemm in the getrs-trsm iterations. This doesn't directly introduce errors into the results but may increase the number of iterations needed for convergence. The effectiveness of low-precision quantization methods can be established if the additional time introduced by low-precision int quantization is less than that of using fp32/fp16 when considering the time added by the increased number of iterations until convergence.

\subsection{Advantages}
Based on our research, the current low-precision quantization solutions are limited to using a single computation method throughout, lacking the ability to balance errors and performance. This paper introduces sparsity for the first time in low-precision dense matrix computations to reduce computational costs.

Besides, the core residual compensation method is carried out after quantization, except for the row-by-column quantization method used in this article. It is also convenient to adopt quantization strategy according to different shapes of input matrix to obtain lower calculation error

The algorithm is not only applicable to Nvidia-GPU but is also easily extendable to other computing devices. It performs well on CPUs and edge computing devices. Additionally, due to the error threshold, the algorithm's computational cost is related to the chosen error threshold. An appropriate threshold within the quantization error range can be selected to meet different error requirements. Moreover, the algorithm automatically chooses the faster execution between sparse matrix multiplication and dense matrix multiplication, avoiding unnecessary computational overhead.

\subsection{Limitation}

Sparse matrices often exhibit significant acceleration only at higher sparsity levels. In some cases, the sparsity obtained through error threshold reduction may not be sufficient to invoke sparse matrix multiplication. As a result, the algorithm continues to execute the complete residual compensation algorithm. Additionally, due to the constraints of low precision, there are still some challenges in applying the algorithm directly to numerical computing tasks.

\section{Related works}
\textbf{In machine learning}:Jung et al. \cite{ref40} proposed quantized pre-training networks in the field of machine learning, while Y. Chen et al.  \cite{ref39} employed residual vector quantization to accelerate the nearest neighbor search process. Lee et al.  \cite{ref41} utilized residual quantization in AR modeling for image processing. Guang Li et al.  \cite{ref31} introduced the use of multiple quantized residuals for error compensation, achieving a noticeable improvement in accuracy at a considerable computational cost (1 to 11 times that of direct quantization). Boyuan Feng et al.  \cite{ref32} implemented arbitrary-length low-precision matrix computations using int1 matrix bitwise operations supported by the Nvidia Ampere architecture, providing support for convolutional neural networks. 

\textbf{In traditional numerical computing}:low-precision calculations have found widespread application. A. Haidar et al. \cite{ref11,ref33} conducted a series of investigations into low-precision mathematical functions, demonstrating the feasibility of mixed precision in high-performance numerical computing. Li, Xiaoye S, Abdelfattah et al. \cite{ref34,ref35} conducted research on low-precision matrix mathematical libraries. Notably, Y. Saad et al. introduced an iterative method using GMRES for solving linear systems, minimizing the norm of the residual vector on the Krylov subspace at each step \cite{ref36,ref37}. Additionally, Zhao, Yuwen et al. \cite{ref38}explored a novel mixed-precision optimization technique, employing a "high-precision computation, low-precision communication" strategy, in the context of applying mixed precision to the fast Fourier transform.

\section{Conclusion}

This paper introduces an error improvement algorithm for quantizing matrix multiplication under low-precision integer data. The algorithm leverages the additional computational efficiency gained from optimizing residual error compensation using sparse matrices generated by an error threshold, and calculation accuracy is improved by combining vector-wise method. Experimental results demonstrate significant acceleration and precision improvement under specific error thresholds, showcasing promising applications in both traditional machine learning and numerical computing domains. Future research will explore the acceleration effects of generating various sparse matrix formats and automatic tuning techniques for different applications based on error threshold values.


\end{document}